%%%%%%%%%%%%%%%%%%%%%%%%%%%%%%%%%%%%%%%%%%%%%%%%%%%%%%%%%%%%%%%%%%%%%%%%%%%
%% Gelfand, I. M.; MacPherson, R. D.
%% 
%% A combinatorial formula for the Pontrjagin classes
%% 
%% A combinatorial formula for the Pontrjagin classes of a triangulated 
%%   manifold is given. The main ingredients are oriented matroid theory 
%%   and a modified formulation of Chern-Weil theory.
%% 
%% publ:  Bull. Amer. Math. Soc. (N.S.) 26(1992) no. 2
%% pp:    304-309
%% type:  Research Announcement        markup: amstex    file size: 23K
%% 
%% copyright: American Math. Society copyright; see end of article
%% 
%% Include files necessary for this article: bull-ppt.tex
%% 
%%%%%%%%%%%%%%%%%%%%%%%%%%%%%%%%%%%%%%%%%%%%%%%%%%%%%%%%%%%%%%%%%%%%%%%%%%%
\input amstex 
\documentstyle{amsppt}
\input bull-ppt
\keyedby{bull282/lbd}
\define\1{\overline}
\define\2{\partial}

\define\al{\alpha}

\define\vare{\varepsilon}

\define\si{\sigma}

\define\Om{\Omega}

\topmatter
\cvol{26}
\cvolyear{1992}
\cmonth{April}
\cyear{1992}
\cvolno{2}
\cpgs{304-309}
\title A combinatorial formula for the Pontrjagin 
classes\endtitle
\author I. M. Gelfand and R. D. MacPherson\endauthor
\address Department of Mathematics, Rutgers University, 
New Brunswick,
New Jersey 08903\endaddress
\address Department of Mathematics, Massachusetts 
Institute of
Technology, Cambridge, Massachusetts 02139\endaddress
\date January 30, 1990 and, in revised form, August 20, 
1991\enddate
\abstract A combinatorial formula for the Pontrjagin 
classes of a triangulated
manifold is given. The main ingredients are oriented 
matroid theory
and a modified formulation of Chern-Weil theory.\endabstract
\endtopmatter

\document
\heading 1. Introduction\endheading
\par
The problem of finding a combinatorial formula for the 
Pontrjagin classes of a
polyhedral manifold $X$ is forty-five years old and has 
stimulated much
research (see \cite{M} and \cite{L} for references). 
Chern-Weil theory provides
a formula for the Pontrjagin classes of a Riemannian 
manifold: They are
represented by differential forms that measure certain 
types of curvature of the
manifold. The problem is to find an analogous theory for 
polyhedra, which have
``infinite curvature'' at the corners.
\par
In this note, we announce a formula that holds in all 
dimensions, is completely
explicit, and can be calculated using combinatorial 
constructions and the
operations of finite-dimensional linear algebra over $\Bbb 
Q$. For each $i$,
the formula gives a rational simplicial cycle $\zeta_i$ in 
the barycentric
subdivision of $X$, whose Poincar\'e dual represents the 
$i\roman{th}$ inverse
Pontrjagin class $\tilde p_i(X)$. (The inverse Pontrjagin 
classes of $X$ are
defined from the usual Pontrjagin classes $p_i(X)$ by
$(1+p_1(X)+p_2(X)+\dotsb)\smile(1+\tilde p_1(X)+\tilde 
p_2(X)+\dotsb)=1$. Like
the $p_i$, the $\tilde p_i$ generate the Pontrjagin ring.) 
The cycle $\zeta_i$
depends on the choice of certain additional combinatorial 
data called a {\it
fixing cycle\/}. We think of a fixing cycle as a 
combinatorial analogue of a
smooth structure on $X$. In fact, a smooth structure on 
$X$ induces a canonical
fixing cycle.
\par
Two ideas make this formula possible: the systematic 
exploitation of oriented
matroids and a reformulation of Chern-Weil curvature theory.
\par
The only other general and explicit combinatorial formula 
for the Pontrjagin
classes is Cheeger's \cite{C}. It uses the asymptotics of 
the spectrum of a
differential operator, so it is difficult to compute and 
its rationality
properties are not clear. However, Cheeger's formula is 
clearly the best one
for the context of the Hodge operator constructed from the 
metric.
\heading 2. Review of oriented matroids\endheading
\par
Oriented matroid theory is a well-developed branch of 
combinatorics with many
applications. One of our goals is to establish a link 
between this theory and
differential topology. We include here only the 
definitions from oriented
matroid theory that we need. The best general reference is 
\cite{BLSWZ}.
\dfn{Definition} An {\it oriented matroid\/} $M$ is a 
finite set $V$ called the
{\it elements\/} of $M$ together with a finite collection 
of maps $c_i\colon
V\to\{-,0,+\}$ called the {\it covectors\/} of $M$ subject 
to the following
axioms:
\roster
\item "1." The constant function with value 0 is a covector.
\item "2." If $c$ is a covector, then $-c$ is a covector 
where $-(-)=+$,
$-(0)=0$, and $-(+)=-$.
\item "3." If $c$ and $d$ are covectors, then $c\circ d$ 
is a covector where
$c\circ d$ is defined by
$$c\circ d(v)=\cases
c(v)&\roman{if}\ c(v)\not=0,\\
d(v)&\roman{otherwise}.\endcases$$
\item "4." For all covectors $c$ and $d$, if $v$ is an 
element such that $c(v)=+$
and $d(v)=-$, then there exists a covector $e$ such that
\item "$\phantom{4.}$" $\bullet$ $e(v)=0$;
\item "$\phantom{4.}$" $\bullet$ If $c(w)=d(w)=0$, then 
$e(w)=0;$
\item "$\phantom{4.}$" $\bullet$ If $c(w)\not=-$, 
$d(w)\not=-$ but $c(w)$ and $d(w)$ are
not both 0, then$\phantom{\bullet\ }\ e(w)=+$;
\item "$\phantom{4.}$" $\bullet$ If $c(w)\not=+$, 
$d(w)\not=+$ but $c(w)$ and $d(w)$ are
not both 0, then$\phantom{\bullet\ }\ e(w)=-$.
\endroster
\enddfn
\par
The idea behind this definition is that an oriented 
matroid is a combinatorial
abstraction of a finite set $V$ of vectors (which are not 
assumed to be
distinct or to be nonzero) in a finite-dimensional real 
vector space $W$. The
covectors $c$ correspond to linear functionals $\1c\colon 
W\to\Bbb R$, but
$c(v)$ remembers only whether $\1c(v)$ is negative, zero, 
or positive.
\par
An element $v$ of $V$ is {\it nonzero\/} if $c(v)\not=0$ 
for some covector $c$.
A subset $\{v_1,v_2,\dotsc,v_j\}$ of $V$ is said to be 
{\it independent\/} if
there exists a set of covectors $\{c_1,c_2,\dotsc,c_j\}$ 
such that $c_i(v_k)=0$
if and only if $i\not=k$. The {\it rank\/} of $x$ is the 
cardinality of any
(and hence every) maximal independent subset of $V$. The 
{\it convex hull\/} of
a set $S$ of elements is $\{v\in V|-\in c(S)$ if 
$c(v)=-\}$. Suppose $N$ and
$M$ are two matroids with the same set of elements $V$. We 
say that $N$ is a
{\it strong quotient\/} of $M$, symbolized $M\Rightarrow 
N$, if every covector
of $N$ is a covector of $M$. If $M$ and $N$ have the same 
rank, we say that the
matroid $N$ is a {\it weak specialization\/} of $M$, 
symbolized
$M\rightsquigarrow N$ if every covector of $N$ is obtained 
from some covector
of $M$ by setting nonzero values equal to zero.
\heading 3. The formula\endheading
\par
Let $X$ be a simplicial manifold of dimension $n$. For 
simplicity, we assume
that $X$ is oriented (with orientation class $[X])$ and 
that $n$ is odd. The
modifications necessary for the general case are noted at 
the end.
\dfn{Definition} The {\it associated complex\/} $Z$ of $X$ 
is the simplicial
complex constructed as follows: The vertices of $Z$ are 
quadruples
$(\Delta,t,y,z)$ where
\roster
\item "$\bullet$" $\Delta\subset V$ is a simplex of $X$, 
where $V$ is the set
of vertices of $X$.
\item "$\bullet$" $t,y$, and $z$ are oriented matroids of 
rank $n+1,2$, and 1
whose set of elements is $V$.
\item "$\bullet$" The matroid $t$ has a covector that does 
not take the
value$-$on any element of $V$.
\item "$\bullet$" The simplex $\Delta$ is related to the 
matroid $t$ by the
following two conditions:
\item "\phantom{$\bullet$}" 1. The nonzero elements of $t$ 
are exactly the
vertices of the star $\roman{St}\,\Delta$ of $\Delta$.
\item "\phantom{$\bullet$}" 2. For each simplex $\Delta'$ in
$\roman{St}\,\Delta$, let $V(\Delta')$ be the set of 
vertices of $\Delta'$.
Then $V(\Delta')$ is linearly independent in $t$ and the 
set of nonzero
elements of $t$ in the convex hull of $V(\Delta')$ is just 
$V(\Delta')$ itself.
\item "$\bullet$" We have strong quotients $t\Rightarrow 
y\Rightarrow z$.
\endroster
The $k$-simplices of $Z$ are diagrams of weak 
specializations and inclusions
$$\matrix
t_0&\rightsquigarrow&t_1&\rightsquigarrow&\dotsb&%
\rightsquigarrow&t_k\\
\Downarrow&&\Downarrow&&&&\Downarrow\\
y_0&\rightsquigarrow&y_1&\rightsquigarrow&\dotsb&%
\rightsquigarrow&y_k\\
\Downarrow&&\Downarrow&&&&\Downarrow\\
z_0&\rightsquigarrow&z_1&\rightsquigarrow&\dotsb&%
\rightsquigarrow&z_k\\\noalign{\vskip 6pt}
\Delta_0&\subseteq&\Delta_1&\subseteq&\dotsb&\subseteq&%
\Delta_k\endmatrix$$
\enddfn
\par
If we delete the matroids $z$ {\it resp.\/} all matroids 
$(t,y$, and $z)$ in
this definition, we get additional associated simplicial 
complexes, which we
denote by $Y$ {\it resp.\/} $\widetilde X$, equipped with 
simplicial maps
$Z\overset{\rho}\to{\to}Y\overset{\pi}\to{\to}\widetilde 
X$. Note that
$\widetilde X$ is just the barycentric subdivision of $X$.
\rem{Remarks} The matroids $t$ are combinatorial 
abstractions of
$T_pX\oplus\Bbb R$ where $T_pX$ is the tangent space to 
$X$ at a point $p$ in
the simplex $\Delta$. The map 
$Y\overset{\pi}\to{\to}\widetilde X$ is an
analogue of the Grassmannian bundle of two-dimensional 
quotients of
$TX\oplus\Bbb R$, and $Z\overset{\rho}\to{\to}Y$ is the 
circle bundle of one
dimensional quotients of the two plane.
\endrem
\thm{Proposition 1} The map $\rho\colon Z\to Y$ is 
topologically a fibration
with a circle as fiber.
\ethm
\demo{Idea of proof} That the fibers over the vertices of 
$Y$ are circles is a
special case of the Folkman-Lawrence representation 
theorem for the oriented
matroid $y$.
\par
We now give a combinatorial formula for the first Chern 
class of a triangulated
circle bundle. Let $\scr O$ be the local system on $Y$ 
with fiber $\Bbb Q$ and
twisting given by the fiber orientation of $Z$. Define a 
1-cocycle $\Theta$ on
$Z$ with coefficients in $\rho^\star \scr O$ as follows: 
For each vertex $v$ of $Y$,
$\Theta|\rho^{-1}v$ is the class that integrates to 1 
around the circle and has
the same value on each 1-simplex. Having fixed this, for 
each edge $e$ of $Y$,
$\Theta|\rho^{-1}e$ is the cocycle such that the sum of 
the squares of the
coefficients is minimum. It is rational, since the problem 
of minimizing a
quadratic expression subject to linear constraints (the 
cocycle condition) can
be solved by linear equations. Now, define the 2-cocycle 
$\Omega$ on $Y$ with
coefficients in $\scr O$ by $\rho^\star\Omega=\delta\Theta$.
\enddemo
\thm{Proposition 2} The cohomology class $\{\Omega\}$ is 
the first Chern class
of the circle bundle $Z$.
\ethm
\demo{Proof} By construction, $\{\Omega\}$ is the 
$\delta_2$ differential
(transgression) of the fiber orientation in the spectral 
sequence for $Z$.
\enddemo
\dfn{Definition} A {\it fixing cycle\/} for $X$ is a 
$(3n-2)$-cycle $\phi\in
Z_{3n-2}(Y,\Bbb Z)$ such that 
$\pi_\star(\Om^{n-1}\frown\phi)=[X]$.
\enddfn
\rem{Remark} A fixing cycle is the combinatorial analogue 
of an orientation
class of the Grassmannian bundle. Unfortunately, it is not 
unique, as its
differential analogue is. The idea evolved from the 
configuration data of [GGL]
and [M].
\endrem
\thm{Theorem 1} Let $\phi$ be a fixing cycle for $X$. Then
$$\tilde p_i(X)\frown[X]=(-1)^i\pi_\star((\tfrac12\Om)^{n+
2i-1}\frown\phi).$$
\ethm
\par
Note that since $n$ is odd, $\Om^{n+2i-1}$ is a cocycle 
with coefficients in
$\Bbb Q$ since $O^{n+2i-1}=\Bbb Q$.
\rem{Remark} This is a cycle level formula, since the 
operations in simplicial
(co)ho-\linebreak mology involved in the right-hand side 
$(\pi_\star,\smile,\frown)$, and
in the definitions of $\Omega$, $\rho^\star$ and $\delta$ 
are all chain level
operations. The complexity of the formula is contained in 
the formulas for
these operations and the construction of the simplicial 
complexes $Z$ and $Y$.
Given the fixing cycle $\phi$, it is a purely local 
formula: the value in an
open set $U$ of $\widetilde X$ depends only on the 
combinatorial structure of
$X$ inside $U$ and on $\phi|\pi^{-1}U$.
\endrem
\heading 4. Construction of the fixing cycle\endheading
\par
First, we study the structure of the auxiliary complex 
$Y$. For any simplex
$\Delta$ in $X$, let $U_\Delta$ be the set of diagrams of 
oriented matroids
$y\Rightarrow t$ with $t$ satisfying conditions 1 and 2 of 
Definition 1 with
respect to $\Delta$. This is a poset by the specialization 
$(\rightsquigarrow)$
ordering on the oriented matroids in $U_\Delta$. We denote 
the order complex of
a poset $P$ by $\roman{Cx}\,P$. The open dual cell of 
$\Delta$ in $\widetilde
X$ is denoted $D\Delta$ (so $\widetilde X$ is the disjoint 
union of the
$D\Delta)$.
\thm{Proposition 3} The subcomplex of $Y$ lying over 
$D\Delta$ is canonically
homeomorphic to $\roman{Cx}\,U_\Delta\times D\Delta$.
\ethm
\par
We denote the homeomorphism of the proposition by
$c^\Delta\colon\roman{Cx}\,U_\Delta\times D\Delta\to Y$.
\rem{Remark} Suppose $\Delta\subset\Delta'$. Then the edge 
of
$c(\roman{C}\times U_\Delta\times D\Delta)$ is glued to
$c(\roman{C}\times U_{\Delta'}\times D\Delta')$ by the map
$\roman{C}\times U_\Delta\to\roman{C}\times U_{\Delta'}$ 
induced by the map of posets
$U_\Delta\to U_{\Delta'}$ defined on the matroid level by 
setting all elements
of $V$ in $\roman{St}\,\Delta$ but not in 
$\roman{St}\,\Delta'$ equal to zero.
\par
A {\it smooth structure\/} on $X$ is a homeomorphism 
$\si\colon X\to M$ to a
smooth manifold that is differentiable on each closed 
simplex in $X$. Define
$\scr Y$ to be the Grassmannian bundle whose fiber over 
$x\in M$ is the space of
$n-1$-dimensional subspaces $F^{n-1}\subset T_xM\oplus\Bbb 
R$. The map
$Y\overset{\pi}\to{\to}\widetilde X$ is a sort of a 
``combinatorial model'' for
the map $\scr Y\to M$.
\par
Let $\scr Y_\Delta$ be the part of $\scr Y$ lying over 
$\si(\Delta)$ for a
simplex $\Delta\subset X$. Any point $y$ in $\scr 
Y_\Delta$ determines an
element of $U_\Delta$ as follows: Let $\si(x)$ be the 
image of $y$ in $M$.
There is a unique embedding 
$e\colon\roman{St}\,\Delta\hookrightarrow T_xM$
that is linear on each simplex, takes $x$ to 0, and 
satisfies
$d(e|\Delta')=d(\si|\Delta')$ for each simplex $\Delta'$ in
$\roman{St}\,\Delta$. Now map $V$ into $T_xM\oplus\Bbb R$ 
by using the
embedding 
$\roman{St}\,\Delta\overset{e}\to{\to}T_xM\overset{\times
1}\to{\to}T_xM\oplus\Bbb R$ for vertices in 
$\roman{St}\,\Delta$ and mapping
all other vertices to zero. This gives a representation of 
the oriented matroid
$t$. The oriented matroid $y$ is represented by projecting 
the images of these
vertices into $(T_xM\oplus\Bbb R)/F^{n-1}$. By this 
construction, $\scr
Y_\Delta$ is decomposed into pieces indexed by elements of 
$U_\Delta$. One can
see from stratified transversality theory that if $\si$ is 
generic, this
decomposition can be refined to a Whitney stratification 
$\scr
Y_\Delta=\bigcup_\al\,S_\al$ that is transverse to the 
boundary. By
construction, each stratum $S_\al$ determines an element 
$u(S_\al)\in U_\Delta$
by which piece of $\scr Y_\Delta$ it lies in.
\par
A {\it full flag of strata\/} in a manifold is a set 
$S=S_0,S_1,\dotsc,S_d$
where the closure of $S_i$ contains $S_{i-1}$, the 
dimension of $S_i$ is $i$,
and $d$ is the dimension of the manifold. If the manifold 
is oriented, the {\it
sign\/} $\vare S$ of $S$ is defined as follows: Map a 
$d$-simplex with vertices
$v_0,\dotsc,v_d$ into the manifold so that the vertex 
$v_0$ goes to $S_0$, the
edge $v_0v_1$ goes to $S_1$, and so on. Then $\vare S=+1$ 
if the orientation of
the simplex agrees with the orientation of the manifold, 
and $\vare S=-1$
otherwise. If $S$ is a full flag of the strata in $\scr 
Y_\Delta$, then denote
by $u(S)$ the oriented simplex in $\roman{C}\times 
U_\delta$ with vertices
$u(S_0),u(S_1),\dotsc$.
\par
For each simplex $\Delta$ of $X$, choose orientations 
$[\Delta]$ of $\Delta$
and $[D\Delta]$ of $D\Delta$ whose cross product is the 
orientation of $X$.
Orient $\scr Y_\Delta$ by the cross product of $[\Delta]$ 
and the standard
orientation of the Grassmannian of $(n-1)$-planes in $(n+
1)$ space (remember
that $n$ is odd).
\endrem
\thm{Theorem 2} The generic smoothing $\si$ induces a 
fixing cycle $\phi$ by
the formula
$$\phi=\sum_\Delta\sum_S\vare(S)c_\star^\Delta(u(\delta)%
\times[D\Delta])$$
where the first sum is over all simplices $\Delta$ of $X$ 
and the second over
all full flags of strata in $\scr Y_\Delta$.
\ethm
\par
The idea of the proof is to construct a continuous map 
$f\colon\scr Y\to Y$ so
that $\phi=f_\star[\scr Y]$.
\heading 5. An alternative form of Chern-Weil 
theory\endheading
\par
Let $E$ be a vector bundle with a connection over a 
differentiable manifold
$M$. Chern-Weil theory gives a formula for the Pontrjagin 
classes of $M$ as a
sum of terms, each of which is a product of curvature 
2-forms $\Om$ multiplied
by a pattern reflecting the structure of the Lie algebra 
of $\roman{Gl}(n,\Bbb
R)$. Finding a combinatorial analogue of $\Om$ is 
possible, but it is a
singular current. The difficulty in finding a 
combinatorial analogue for
Chern-Weil theory is regularizing the products.
\par
The combinatorial formula of this paper is an analogue of 
another form of
Chern-Weil theory, which we now describe. Let $e$ be the 
fiber dimension of $E$
and assume that it is even. Let $\pi\colon\scr Y\to M$ be 
the Grassmannian
bundle of $(e-2)$-planes in $E$ and $\rho\colon\scr 
Z\to\scr Y$ be the
principle circle bundle of the tautological quotient 
2-plane bundle $\xi$ over
$\scr Y$. The connection on $E$ induces a one form 
$\Theta$ on $\scr Z$ with
coefficients twisted by the orientation sheaf $\scr O$ of 
$\scr Z$ and a
curvature form $\Omega$ on $\scr Y$ defined by 
$\rho^\star\Om=d\Theta$.
\thm{Proposition 4} $\tilde 
p_i(E)=(-1)^i\pi_\star\Om^{(e-2+2i)}$ where
$\pi_\star$ represents integration over the fiber.
\ethm
\par
Proposition 4 is proved by an algebraic manipulation using 
only the Whitney sum
formula applied to $\pi^\star E=\xi\oplus\xi^\perp$, the 
vanishing of high
Pontrjagin classes of a low-dimensional bundle, and the 
projection formula.
Theorem 1 is a combinatorial analogue of this formula, 
where $E$ is $TM\oplus
1$.
\subheading{Orientations and dimensions} Suppose that $X$ 
is not orientable
and/or not odd dimensional. Let $\scr D$ be the 
orientation local system of
$X$, so $[X]\in H_n(X,\scr D)$. The fixing cycle should 
lie in homology with
twisted coefficients: $\phi\in Z_{3n-2}(Y,\pi^\star\scr 
D\otimes\scr
O^{\otimes(n-1)})$. The construction of $\phi$ in \S2 
still works because $\scr
Y$ has orientation sheaf $\pi^\star\scr D\otimes\scr 
O^{\otimes(n-1)}$ where
$\scr O$ is the orientation sheaf $\scr Z$.


\Refs\ra\key{BLSWZ}

\ref\key{BLSWZ} 
\by A. Bj\"orner, M. Las Vergnas, B. Sturmfels, N. White, 
and G.
Ziegler \paper Oriented matroids 
\inbook Encyclopedia Math. Appl.
\publ Cambridge Univ. Press
\yr 1992 \endref

\ref\key{C} 
\by J. Cheeger \paper Spectral geometry of singular 
Riemannian spaces 
\jour J. Differential Geom. \vol 18 
\yr 1983 
\pages 575--657 \endref

\ref\key{GGL} 
\by I. Gabrielov, I. Gelfand, and M. Losik \paper 
Combinatorial
calculation of characteristic classes  
\jour Funktsional Anal. i Prilozhen. \vol 9 
\yr 1975 
\pages 54--55
\moreref
\pages no. 2, 12--28
\moreref 
\pages no. 3, 5--26 \endref

\ref\key{M} 
\by R. MacPherson \paper The combinatorial formula of 
Garielov, Gelfand,
and Loskik for the first Pontrjagin class  
\paperinfo S\'eminaire Bourbaki No. 497
\inbook Lecture Notes in Math. vol. 667
\nofrills\publ Springer
\publaddr Heidelburg
\yr 1977 \endref

\ref\key{N} 
\by N. Levitt \paper Grassmannians and gauss maps in 
piecewise-linear
topology  
\inbook Lecture Notes in Math.
\vol 1366
\publ Springer
\publaddr Heidelburg
\yr 1989 \endref
\endRefs
\enddocument